\widowpenalty=15000
\overfullrule=0pt
\pretolerance=-1
\tolerance=2500
\doublehyphendemerits 50000
\finalhyphendemerits 25000
\adjdemerits 50000
\hbadness 1500
\abovedisplayskip=7pt plus 5pt minus 6pt
\belowdisplayskip=7pt plus 5pt minus 6pt
\hsize39pc \vsize52pc
\baselineskip11.5pt
\magnification=1200
\parindent 0pt
\def\exp{{\rm \hskip1.5pt  exp \hskip1.5pt }}

\font\gross=cmbx12
\font\grossrm=cmr12
\font\gr=cmr17


%
%
\def\parein#1#2{\par\noindent\rlap{\rm {#2}}\parindent#1\hang\indent\ignorespaces}
\def\paraus{\parindent 0pt\par}

\def\G{\Gamma}

\def\l{\ell}

\def \1{\backslash}
%
%
%
\output={
\ifodd\pageno\hoffset=0pt
\else
\hoffset=0truecm\fi\plainoutput}
\font\authorfont=cmcsc10 at 12pt
\font\titlefont=ptmb at 13pt
\font\affiliationfont=cmti12 at 11pt
\font\headingfont=cmr12 
\font\gross=cmbx12

\font\ninecaps=cmcsc9
\font\ninebf=cmbx9
\font\nineit=cmti9
\font\ninerm=cmr9
\font\eightrm=cmr8
\font\eightit=cmti8

\def\re{{\rm Re}\,}

\def\la{{\bf L}}

\def\chyp #1 { C_{{\rm hyp}, #1 }}

\abovedisplayskip=8pt plus 5pt minus 3pt
\belowdisplayskip=8pt plus 5pt minus 3pt
 at 11truept

\def\re{{\rm Re}\,}

\def\sgn{{\rm sgn}\,}
\def\fh{{\bf F}}
\def\G{\Gamma}

\def\erfc{{\rm Erfc}\,}
\def\erf{{\rm Erf}\,}
%
%
%
%
%
%
\output={
\ifodd\pageno\hoffset=0pt
\else
\hoffset=0truecm\fi\plainoutput}
\font\authorfont=cmcsc10 at 12pt
\font\titlefont=ptmb at 13pt
\font\affiliationfont=cmti12 at 11pt
\font\headingfont=cmbx12
\font\paragraphfont=cmbx10
\font\rfont=cmbx10
\font\gross=cmbx12

\font\ninecaps=cmcsc9
\font\ninebf=cmbx9
\font\nineit=cmti9
\font\ninerm=cmr9
\font\eightrm=cmr8
\font\eightit=cmti8

\def\re{{\rm Re}\,}

\abovedisplayskip=8pt plus 5pt minus 3pt
\belowdisplayskip=8pt plus 5pt minus 3pt
\baselineskip=15pt
%
\centerline{\titlefont ON THE VALUATION OF PARIS OPTIONS:}
\centerline{\titlefont FOUNDATIONAL RESULTS}
\baselineskip=11.5pt
\vskip.3cm
\centerline{\authorfont Michael  Schr\"oder}
\vskip.2cm
\centerline{\affiliationfont Lehrstuhl Mathematik III}
\centerline{\affiliationfont Seminargeb\"aude A5,
Universit\"at Mannheim,  D--68131 Mannheim}
\vskip.4cm 
\centerline{
\vbox{
\hsize=10.8cm
\baselineskip=9.3pt
\ninerm 
This paper addresses the valuation of the Paris barrier 
options of [{\ninebf CJY}] using the Laplace transform approach. The 
notion of Paris options is extended such that their valuation 
is possible at any point during their lifespan.
The Laplace transforms of [{\ninebf CJY}] are modified when 
necessary, and their basic analytic properties are discussed.
}
}
\vskip.4cm
{\paragraphfont 1.\quad Introduction:}\quad This paper is the first part of 
a report on the valuation of the Paris form of barrier options proposed 
in [{\bf CJY}]. In their standard form European--style barrier options 
come as puts or calls that are activated or deactivated as soon as 
their underlying hits a pre\-speci\-fied barrier level. The new 
idea of the Paris barrier options for cushioning this abruptness 
is to introduce a systematic delay for these consequences of hitting 
the barrier to become effective. 
\medskip
To fix ideas consider the  Paris down--and--in call option. Suppose the 
underlying is moving downwards and hits the barrier level. This Paris 
option  is then not immediately activated. Instead, a clock starts 
ticking at the very moment the underlying passes below the barrier level. 
It will tick for $D$ consecutive periods in time. Only then, with this 
delay of $D$ periods the option will be activated if the following two
conditions hold. First, one is still during the life span of the option. 
Second, during any of these last $D$ periods of time the underlying has 
been staying below the barrier level. 
To further illustrate this Paris option, suppose $D$ periods before its 
maturity the option is not yet activated and no clock has started ticking. 
Then  the option will not be activated during its life time. It is 
already worthless. 
\medskip
This paper addresses foundational questions in valuating 
Paris options. For instance, a principal difficulty  had been pointed 
out to the author by Pliska.  Suppose the clock is already ticking for
a Paris option, how should one valuate this option? The necessary 
extensions of [{\bf CJY}] are discussed in Part~I. For valuating 
Paris options the Laplace transform approach is adopted. Indeed, Laplace 
transforms of their prices have been computed in [{\bf CJY}], and the 
discussion of Part~II is based on their results. However, pursuing 
Pliska's suggestions we found it necessary to modify one of these two 
Laplace transforms. So I am very grateful to Pliska for his remark. 
In this conncetion I also wish to thank Yor for making  available to 
me corrections to [{\bf CJY}] collected by Jeanblanc--Picqu\'e. 
\medskip
This status quo ante set by [{\bf CJY}] is extended in Part~III where 
a number of insights into the analytical structure of the above Laplace 
transforms are discussed. More precisely, these Laplace transforms are 
constructed using a function $\Psi$ whose analytical properties are 
analyzed. In particular the series of this function is derived and also 
an asymptotic expansion for it is proved that is uniform on the 
right--half plane. These results are fundamental for the Laplace 
inversion of the above Laplace transforms.
\medskip
Recall that the Laplace transform generally consists of two stages. 
As a first step, a 
suitably nice function on the positive real line is transformed 
into a complex--valued function. This Laplace transform exists on a 
half--plane sufficiently deep within the right complex half--plane  
and defines there a holomorphic, i.e., complex analytic, function.
As a second step, the Laplace transform has to be inverted to give 
the desired function on the positive real line. The natural methods 
in this step are complex analytic. Relevant concepts for the
endeavours ahead are recalled in Part~0.
\medskip
To proceed, recall how for valuating barrier options the 
relative position of strike price and barrier play a role. 
For Paris options this list has to be extended 
and we propose the following notions. If the barrier is equal to 
or below the strike price we have a standard situation, otherwise
we have a perverse situation. Also the relative position of today's 
price of the underlying to the barrier level matters. We have the 
first situation if today's price of the underlying is equal to or 
above the barrier level. 
If today's price of the underlying is below the barrier level, so
that the clock for Paris option has started ticking, we have the 
second situation. These situations combine into four possible cases. 
In each of these cases there is a different Laplace transform to be 
inverted. 
\medskip
With this program in mind the reader now may wish not only to take 
a look at this paper but also at its companion [{\bf SP}]. The focus
of this last paper is on Laplace inversion. It addresses valuating 
the Paris options in the above first standard case. However, its 
results about Laplace inversion are typical and 
fundamental for the remaining three cases.  
\nopagenumbers
\def\Verfasser{{\ninecaps Michael Schr\"oder}}
\def\Kapiteltitel{{\ninecaps Valuation of Paris options}}
\headline={\vbox{\line{\ninebf \ifodd\pageno\hss \Kapiteltitel\hss
\folio
\else
\folio\hss\Verfasser\hss\fi}
}\hss}
\bigskip
\vbox{\ninerm
\baselineskip=9.4pt
{\ninebf Acknowlegements:}\quad In addition to Professors Pliska and
Yor already mentioned in the introduction, I would like to thank 
in Mannheim Professor R. Weissauer, Dr.~U. Weselmann, and Dr.~D. Fulea 
for their interest in and support of this project.}
\goodbreak
\bigskip
\centerline{\headingfont Part 0\qquad  Epitome of methods}
\vskip.3cm
{\paragraphfont 2.\quad  Certain classes of  higher transcendental 
functions:}\quad This paragraph recalls pertinent facts about the 
the gamma function and the error function from [{\bf L}], [{\bf D}]. 
\goodbreak
\bigskip
{\paragraphfont Gamma function:}\quad The {\it gamma function\/} $\G$ is
defined on the right complex half--plane by the improper integral:
$$ \G(z)=\int_0^\infty e^{-t} t^{z-1}\, dt\, ,$$
for any complex number $z$ with positive real part. 
It interpolates the factorial in the sense that $\G (n+\!1)=n!$, for any 
non--negative integer $n$, and more generally satisfies the {\it recursion
rule\/} $\G(z\!+\!1)=z\G(z)$. The {\it functional equation\/}
$ \G(z)\cdot \G(1\!-\!z)\cdot\sin(\pi z)=\pi$
in particular shows how it continues to a meromorphic function on the 
whole complex plane with simple poles  in the non--positive integers as its
singularities. It  satisfies the  {\it duplication formula\/}
$2^{2z-1}\G(z) \G( z\!+\!1/2) =\sqrt{\pi}\,\G(2z)$.
\medskip
For any complex numbers $a$, $z$ with positive real part, one has 
the identity: 
$\G(z)=a^z\int_{\raise1.5pt\hbox{$\scriptstyle 0$}}^\infty 
e^{-at}t^{z-1}\, dt$. Indeed, for $a$ real it is obtained by a change 
of variable. Applying the identity theorem the resulting identity is 
then valid on the whole right half--plane.
\goodbreak
\medskip
Unlike the gamma function itself, its reciprocal is  
holomorphic on the whole complex plane. A repesentation as a contour
integral is given by the {\it Hankel formula\/}:
$$ {x^{a-1}\over \G(a)}={1\over 2 \pi i}\int_C e^{xw} w^{-a} \, dw \, ,$$
that is valid for any complex number $x$ with positive real part, and now
for any complex number $a$. Herein $C$ is any of the {\it Hankel contours\/}
$C_{\varepsilon, R}$ with $0<\varepsilon<R$, or their limits for 
$\varepsilon$ converging to zero, defined as the following paths of 
integration: enter from minus infinity 
\input xy
\xyoption{arc}
\xyoption{arrow}
$$
\xy     0;<1mm,0mm>                
                \ar (10,-5);(30,-5),            
                \ar @{-}(30,-5);(50,-5),        
                \ar @{-*\dir{>}}(50,5);(30,5),  
                \ar @{-}(30,5);(10,5),          
                \ar @{.}(50,5);(55,5),	
                \ar 0;(72.06,0),                        
                \ar (55,-12);(55,12),\ar@{},            
                (55,5)*+\dir{*}*+!LU{\scriptstyle i\varepsilon},
                (62.06,0)*+\dir{*}*+!UL{\scriptstyle R},
                (47.94,0)*+\dir{*}*+!UL{\scriptstyle -R},
                (55,0)*\cir<7.06mm>{dr^dl},\ar@{},

\endxy
$$
\vskip-.1cm
\centerline{\eightrm Figure A\enspace The Hankel contour 
$\hbox{\eightit C}_{\scriptstyle \varepsilon,R}$}
\goodbreak
\vskip.4cm
on a parallel to the real axis through $-i\varepsilon$ until 
circle in the orgin with radius $R$ is first hit. Pass counterclockwise on
this circle until the parallel to the real axis through $+i\varepsilon $ is 
hit in the left complex half--plane. Leave on this parallel to minus
infinity. 
\goodbreak
\bigskip
{\paragraphfont Error function:}\quad The {\it complementary error function\/}
$\erfc$ is the analytic function on the complex plane~given~by:
$$ 
\erfc(z)
={2\over \sqrt{\pi\,}}\int_z^\infty e\vbox to 8pt{}^{-w^2}\, dw\, ,
$$
for any complex number $z$, where the path of integration in the last
integral starts at $z$ and moves parallel to the real axis with
the real part of its elements  going to infinity. 
\goodbreak
\bigskip
{\paragraphfont 3.\quad Laplace transform:}\quad This paragraph collects
pertinent facts about the Laplace transform from [{\bf D}].
\bigskip
The class of functions 
considered in the sequel is that of functions of {\it exponential type\/}, 
i.e., of continuous, real--valued functions $f$ on the non--negative real 
line such that there is a real number $a$ for which $\exp(at)f(t)$ is 
bounded for any $t>0$. The {\it Laplace transform\/} is the  linear 
operator $\la$ that associates to any function $f$ of exponential type 
the complex--valued  function  $\la(f)$, analytic on 
a suitable complex half--plane, given by:
$$ \la(f)(z)=\int\nolimits_0^\infty e^{-zt}f(t)\,dt\, ,
$$
for any complex number $z$ with sufficiently big real part.
The {\it convergence theorem\/} is that the Laplace transform 
$\la(f)$ converges on the whole half plane $\{z| \re(z)\ge \re(z_0)\}$
if $\la(f)(z_0)$ exists. So there is a unique 
$-\infty\le \sigma=\sigma_f\le \infty$, the {\it abscissa of 
convergence\/} of $\la(f)$, such that $\la(f)(z)$ converges, for any $z$ with
$\re(z)>\sigma$ and does not exist for $z$ with $\re(z)<\sigma$. 
To fix ideas, consider the function given by $f(t)=\exp(at)$ 
with $a$ any positive real number so that formally 
$\la(f)(z)=
\int_{\raise1.5pt\hbox{$\scriptstyle 0$}}
^{\raise-1.5pt\hbox{$\scriptstyle \infty$}} \exp(-(z\!-\!a)t)\, dt$. 
With the absolute value of $\exp(w)$ being equal to 
$\exp(\re(w))$, for any complex number $w$, this improper integral 
exists iff the real part of $z\!-\!a$ is positive, 
i.e., iff the real part of $z$ is bigger than $a$. In this 
case, $\la(f)(z)$ equals $(z\!-\!a)^{-1}$ and  is analytic on the 
half plane $\{ z| \re (z)>a\}$.
\medskip
That $\la$ is an injection is a consequence of ultimately the 
Weierstra\ss\ approximation theorem. 
Its inverse, the inverse Laplace transform  $\la ^{-1}$, is expressed 
as a contour integral by the {\it complex inversion formula\/} of 
Riemann. For any function $H$ analytic on 
half--planes $\{ \re (z)\!\ge\!x_0\}$ with $x_0$ any sufficiently 
big positive real number, it asserts: 
$$ \la^{-1}(H)(t)={1\over 2\pi i}\int\nolimits_{x_0-i\infty}^{x_0+i\infty}
e^{zt}\cdot H(z)\, dz,$$ 
for any positive real number $t$ if $H$ satisfies a  growth 
condition at infinity such that the above integral exists. In particular
$(\la ^{-1}\circ \la) (f)=f$ if $f$ is of exponential type.  
\medskip
The convolution of two functions $f$ and $g$ of exponential type is the 
function of exponential type given by  
$f*g(t)=\int\!\vbox to 7pt{}_{\! 0}^{\,t} f(u)g(t-u)\, du$ 
for any $t\ge 0$. With $du$ being Lebesgue measure, notice $f*g=g*f$. 
The {\it product theorem\/} for the Laplace transform is that $\la$
transforms convolutions into products:
$$ \la(f*g)(z)=\la(f)(z)\cdot \la(g)(z), $$
for any $z$ with $\re (z)$ sufficiently big. Equivalently, 
$f*g(t)=\la^{-1}( \la(f)\cdot \la(g))(t)$. The proof is a formal
Fubini calculation as soon as one extends functions of exponential type
by zero on the non--positive real line, whence $f*g(t)$ equals  
$\int_0^\infty f(u)g(t-u)\, du$. 
\medskip
There are a number of {\it shifting principles\/} for the Laplace transform.
First, weighting any function $f$ of exponential 
type with the function sending any $t>0$ to $\exp(at)$, where $a$ 
is any real number, introduces a shift by minus $a$ in the Laplace 
transform:
$$ \la(e^{at}f)(z)=\la(f)(z\!-\!a).$$
The effect of weighting the Laplace transform of $f$ by $\exp(az)$ 
depends on the sign of the real number $a$. More precisely, let $a$ 
be any non--negative real number. Using the above complex inversion 
formula, one has:
$$ \la^{-1}\bigl(e^{az}\la(f)(z)\bigr)(t)=f(t+a), $$ 
for any $t>0$. Weighting with $\exp(-az)$ gives the following 
{\it shifting theorem\/}:
$$  \la^{-1}\bigl(e^{-az}\la(f)(z)\bigr)(t)={\bf 1}_{(a,\infty)}(t)f(t-a), $$
for any positive real number $t$. Herein ${\bf 1}_X $ is the 
characteristic function of any set $X$. This is proved calculating the 
Laplace transform of the right hand side. 
\goodbreak 
\medskip
The {\it integration and differentiation principles\/} for the 
Laplace transform describe the behaviour of the Laplace transform
upon multiplication by integral powers of the complex variable $z$. 
Using partial integration in the defining integral of the Laplace 
transform:
$ \la(f^{(n)})(z)
=z^n\cdot \la(f)(z)-  f^{(0)}(0)z^{n-1}\!-\cdots-\! f^{(n-1)}(0)z^{0}$,
for any non--negative integer $n$ and any complex number $z$ in the
half--plane of convergence for $\la(f)$. In particular, 
$$\la\bigl(f^{(n)}\bigr)(z)=z^n\cdot\la(f)(z),$$  
if $f$ is such that for any non--negative integer $k$ less than $n$ its
$k$--th derivative $f^{(k)}$ vanishes at zero. This last formula holds
in general if one considers $n$--fold integration as differentiation of 
the negative order $-n$. More precisely, 
the Laplace inverse
of $z^{-n}\cdot \la(f)(z)$, for any non--negative integer $n$, is given as 
the $n$--fold {\it iterated integral\/} $I_n(f)$ of $f$. This is 
recursively defined by $ I_0(f)=f$, and 
$$ I_{m+1}(f)(t)=\int_0^t I_m(f)(u)\, du\, ,$$
for any positive real number $t$, and any non--negative integer $m$. 
\goodbreak
\vskip.7cm
\centerline{\headingfont Part I\qquad Paris options}
\vskip.3cm
{\paragraphfont 4.\quad Black--Scholes framework for valuating 
contingent claims:}\quad
For our analysis we place ourselves in the Black--Scholes framework and use
the risk--neutral approach to valuating contingent claims. 
In this set--up one has two securities. First there is a riskless security,
a bond, that has the continuously compounding positive interest rate $r$.
Then there is a risky security whose price process $S$ is modelled as 
follows. Start with a complete probability space equipped with the standard 
filtration of a standard Brownian motion  on it that has the time set 
$[0,\infty)$. On this filtered space one has the {\it risk neutral
measure\/}  $Q$, a probability measure equivalent to the given one, and
 a standard $Q$--Brownian motion $B$ such that $S$ is the strong solution 
of the following stochastic differential equation:
$$ dS_t=(r-\delta)\cdot S_t\cdot dt +\sigma \cdot S_t\cdot dB_t, 
\qquad t\in[0,\infty).$$  
Herein the positive constant $\sigma$ is the volatility of $S$. The 
constant $\delta$ depends on the security modelled. For instance, it is 
zero if $S$ is a non--dividend paying stock.
\goodbreak
\bigskip
{\paragraphfont 5.\quad Paris options:}\quad This paragraph
dicusses the Paris barrier options of [{\bf CJY}, \S 2]. Their new idea for 
cushioning the impact of the underlying hitting the barrier is as follows. 
They require the underlying $S$ to spend a minimum time $D>0$ 
above or below their prespecified barrier $L\ge 0$ before the option is 
knocked in or knocked out.
\medskip
The {\it Paris down--and--in call\/} to be considered in the sequel 
is the following European--style contingent claim on $S$ written at time
$t_0$ and with time to maturity $T$. Its payout at $T$ is that of 
a call on $S$ with exercise price $K$:
$$ (S_T-K)^+=\max\bigl\{ S_T-K,0\big\}\, , $$
if $S_t$ is less than $L$ during a connected subperiod of length at least 
$D$ of the monitoring period $[t_0,T]$. Equivalently, there is a point in 
time $a$ such that the interval $I=(a,a+D)$ is contained in the lifespan 
$[t_0, T]$ of the option and $S_t<L$ for any $t$ in $I$. Otherwise the 
call expires worthless.
\medskip
In the sequel, fix any time $t$ during $[t_0,T]$ at which the Paris  
down--and--in call has not yet been knocked in. The time--$t$ value 
$C_{d,i}$ of this call is then a function of $K$, $\tau=T\!-\!t $, 
$ S_0$, $r$, $\delta$, $L$ and $D$. In the sequel, mostly any of these 
arguments is suppressed from the notation. The valuation problem depends 
on an excursion time $H_{L,t}=H_{L,D,t}^-$ defined in the sequel such that 
the Paris down--and--in call is knocked in before its maturity $T$ if and 
only if $H_{L,t}\le T$. Using the arbitrage pricing principle, the 
time--$t$ price $C_{d,i}$ of the Paris down--and--in call is then  
given by the conditional $Q$--expectation:
$$ C_{d,i}=e^{-r\tau} E_t\Big[ \phi(S_T)\cdot {\bf 1}_{\{ H_{L,t}\le T\}}
\Big],$$
with $\phi(S_T)$ being equal to the plain vanilla call payout 
$(S_T\!-\! K)^+$, and denoting the characteristic function of the event 
$H_{L,t }\le T$ by  ${\bf 1}_{\{ H_{L,t }\le T\}}$.  
\medskip
In defining $H_{L,t}$ at time $t$ two constellations are to be 
distinguished. Namely, either today's security price $S_t$ is equal to 
or above  the barrier $L$, or it is below $L$. 
In the first case with $S_t\ge L$, the Paris down--and--in call is 
knocked in before its maturity $T$ iff $S$ remains smaller 
than $L$ for all points in time of an interval of length at least $D$ 
contained in $[t,T]$. Thus let $H_{L,t}$ in this case denote the 
first time $s$ greater than or equal to $t\!+\!D$ such that $S_u<L$  
for any time $u$ in the interval $(s\!-\!D,s)$.
\medskip
In the second case where $S_t<L$, the security price $S$ has already 
been staying below $L$ for some connected period time during
the lifetime of the option. The down--and--in call is then knocked in 
before its time to maturity $T$ if the following happens. The security 
price $S$ continues to stay below $L$ also during the period of time 
from today until the point in time $t\!+\!d$ with $d$ smaller than $D$,  
i.e., one has $S_u<L$ for all $u$ in $(t,t\!+\!d)$. Thus define 
$H_{L,t}=t\!+\!d$ in this case. 
If the level $L$ is hit by $S$ before time $t\!+\!d$, however, the clock 
for the minimum lenght $D$ is restarted. So define $H_{L,t}$ in this case 
to be the first time $s$ greater than or equal to $t\!+\!D$ such that 
$S_u<L$ for any $u$ in the subinterval $(s\!-\!D, \!s)$ of the positive 
real line.  
\medskip
The discussion of Paris barrier options would not be complete without 
having referred to [{\bf CJY}, \S 2] for details on the whole family 
of Paris barrier options. Calls and puts with either of the 
following barrier types: down--and--out, down--and--in, up--and--out, 
up--and--in. 
Valuating these eight types of options is reduced in a standard 
way to the case of the above Paris down--and--in call. Indeed, put--call 
parities are given in [{\bf CJY}, \S 5], and [{\bf CJY}, \S 4.2] reduces 
the valuation of the out--calls to that of the in--calls, mutatis mutandis. 
This paper thus concentrates on valuating the  Paris down--and--in call, 
in the sequel also referred to as {\it the Paris option\/} for simplicity. 
\goodbreak
\vskip.7cm
\centerline{\headingfont Part II \qquad Laplace transforms}   
\vskip.3cm
{\paragraphfont 6.\quad Laplace transforms of the value of the 
Paris option:}\quad This paragraph dis\-cus\-ses the Laplace transforms  
of the densities for the the Paris option of paragraph five. 
\bigskip
{\paragraphfont The basic valuation identity:}\quad This notion of 
density is made precise by by the following {\it basic valuation identity\/}
for the time--$t$ price $C_{d,i}$ of the Paris option:
$$C_{d,i}=e\vbox to9pt{}^{-
        \bigl( r+{\scriptstyle \varpi^2\over\scriptstyle 2}\bigr)
                        \cdot\tau }
\int_{\beta(S_t)}^\infty e^{\varpi x}\cdot \bigl( 
S_t\cdot e^{\sigma x}\!-\!K\bigr) 
\cdot h_b(\tau,x)\, dx\, ,$$
with the {\it Paris option density\/} $h_b$ given by:
$$
h_b(u,y)=
\int_{-\infty}^\infty 
E^*\biggl[ {\bf 1}_{\{ H_{b}^* <u\}}
{1\over \sqrt{ 2\pi\cdot(u\!-\! H_{b}^*)\, }} 
\cdot e\vbox to 9pt{}^{
-{\scriptstyle (x-y)^2\over\scriptstyle 2 \cdot (u\!-\! H_{b}^*)}
              }\, \biggr]d\mu^*(x)\, , $$
for any real numbers $u>0$ and $y$, and where 
$\beta(S_t)=\sigma^{-1}\log(K/S_t)$, where  $H_b^*=H_{L,t}-t$, 
and with $\mu^*=\mu_{b,D,-}$ the measure for Brownian motion 
at time $H^*_b$. 
By construction, $h_b(u,\ )$ is zero for $u$ less than $D$, respectively $d$ 
depending on today's situation. 
\goodbreak
\bigskip
{\paragraphfont Laplace transforms of the densities:}\quad With the 
proof of the basic valuation identity postponed to the end of this 
paragraph, the Laplace
transforms of $h_b$ with respect to time are discussed next. For any
fixed real number $y$ they are recalled to be defined by:
$$ \la\bigl( h_b(\ ,y)\bigr)(z)=\int_0^\infty e^{-zt} h_b(t,y)\, dt \, , $$
for any complex number $z$ with sufficiently big real part. 
Moreover define the~function~$\Psi$~by:
$$ \Psi(z)=\int_0^\infty x\cdot 
e\vbox to 9pt{}^{-{\scriptstyle x^2\over \scriptstyle 2}+zx}dx\, , $$
for any complex number $z$, and choose on the complex plane with the 
non--positive real axis deleted the square root defined using the principal
branch of the logarithm. The computations of [{\bf CJY}, \S\S 5, 8] can 
then be summarized as follows. 
\bigskip
{\rfont Proposition A:}\quad {\it For any real number $y$, the 
Laplace transform of $h_b(\ , y)$ is a holomorphic function on the 
right complex half--plane.}
\bigskip
{\rfont Proposition B:}\quad {\it Suppose 
$b=\sigma^{-1}\log(L/S_t)$ is non--positive. For any fixed real number $y$, 
the Laplace transform of $h_b(\ ,y)$ is given by:     
$$ \la\bigl( h_b(\ ,y)\bigr)(z)=
{e\vbox to 9pt{}^{{\scriptstyle b\over \scriptstyle \sqrt{D\, }}
                \sqrt{ 2Dz\, }
               }
\over
\sqrt{D\,} \sqrt{ 2Dz\, }\Psi(\sqrt{ 2Dz\, })
}
\int_0^\infty 
x\cdot e\vbox to 9pt{}^{-{\scriptstyle x^2\over \scriptstyle 2D} -
                       |b-x-y|\sqrt{ 2z\, }}
dx\, , $$
for any complex number $z$ with positive real part.}
\bigskip
If $b=\sigma^{-1}\log(L/S_t)$ is positive, i.e., when today's security
price $S_t$ is below the barrier $L$, we have modified the Paris option
of [{\bf CJY}]. Today's situation is then such that the price of the 
security has already been staying below the barrier $L$ for a certain 
connected period of time during the lifetime of the Paris options. 
For the Paris option to be knocked in the price thus has to stay 
below the barrier $L$ only for another connected period time of a 
length $d< D$ from today, i.e., $S_u<L$ for all $u$ in $[t, t\!+\!d)$. 
In this way, Paris options now can be valuated at any point of their
monitoring periods. The desirability of this modification
has been  pointed out to me by Pliska. I regard the following 
modification of the valuation results of [{\bf CJY}] as a direct 
result of his suggestions.
\goodbreak
\bigskip
{\rfont Proposition C:}\quad {\it Suppose 
$b=\sigma^{-1}\log(L/S_t)$ is negative. For any real number $y$, 
the Laplace transform of $h_b(\ ,y)$ is given on the right half--plane  
as the following four--term sum of Laplace transforms:}
$$\eqalign{
\la\bigl(h_b(\ ,y)\bigr)=\, & 
\erfc\Bigl( {b\over \sqrt{2d\, }}\Bigr)\cdot \la\Bigl(h_{b,1}(\ , y)\Bigr)
+\la\Bigl(h_{b,2}(\ , y)\Bigr)\cr 
\noalign{\vskip4pt}
&+
\erfc\Bigl( {b\over \sqrt{2d\, }}\Bigr)
\erf\Bigl( {b\over \sqrt{2d\, }}\Bigr)
\cdot \la\Bigl(h_{b,3}f(\ ,y)\Bigr)
+
\erf\Bigl( {b\over \sqrt{2d\, }}\Bigr)
\cdot \la\Bigl(h_{b,4}(\ ,y)\Bigr)\cr 
}
$$
\medskip
Herein $h_{b,k}$ are the functions on the positive real line times 
the real line defined  by:
$$\eqalign{
h_{b,3}(u,y)&={\bf 1}_{(d,\infty)}(u)\cdot {1\over u\!-\!d\!+\!D}
\bigg\{ 
{\sqrt{u\!-\!d\, }\over \sqrt{2\pi\,}}\cdot 
 e\vbox to 12pt{}^{ -{\scriptstyle (y-b)^2\over\scriptstyle 2(u\!-\!d)}  }
\cr 
\noalign{\vskip5pt}
&\hbox to 2cm{}
-{1\over 2}{(y\!-\!b)\sqrt{D}\over \sqrt{u\!-\!d\!+\!D}\,}\cdot  
e\vbox to 12pt{}^{\! -{\scriptstyle (y-b)^2
\over\scriptstyle 2(u\!-\!d\!+\!D)}  }
\erfc\!\biggl( {(y\!-\!b)\sqrt{ D}\over 
\sqrt{2(u\!-\!d)(u\!-\!d\!+\!D)}\,}\biggr)
\bigg\},
\cr \noalign{\vskip5pt}
h_{b,4}(u,y)&={\bf 1}_{(d,\infty)}(u){1\over \sqrt{2\pi u\,}}\cdot
f_{b,u}(y)\qquad \hbox{where}\qquad 
f_{b,u}(y)=e\vbox to 9pt{}^{-{\scriptstyle y^2\over\scriptstyle 2u}  }
- 
e\vbox to 9pt{}^{ -{\scriptstyle (y-2b)^2\over\scriptstyle 2u}  }
,
}$$
for any positive real number $u$ and any real number $y$, while the 
remaining two functions are defined using their Laplace transforms 
as follows:
$$\displaylines{
\la\bigl((h_{b,1}(\ , y)\bigr)(z)
=
\int_0^d 
{e\vbox to 9pt{}^{-zw }\over D\sqrt{2z\, }\Psi(\sqrt{2Dz\, })} 
\biggl( 
\int_0^\infty 
x\cdot e\vbox to 9pt{}^{-{\scriptstyle x^2\over \scriptstyle 2D} -
                       |b-x-y|\sqrt{ 2z\, }}dx
\!\biggr) \mu_b(dw)\, , \cr 
\noalign{\vskip5pt}
\la\bigl((h_{b,2}(\ , y)\bigr)(z)
=
{1\over \sqrt{2\pi d\, }}\int_0^d
{e\vbox to 9pt{}^{-zw }\over \sqrt{2z\, }\Psi(\sqrt{2Dz\, })} 
\biggl( 
\int_{\bf R}e\vbox to 9pt{}^{-|x-y|\sqrt{2z\, }} f_{b,d}(x)
dx
\!\biggr) \mu_b(dw)\, ,\cr 
}
$$
for any positive real number $y$ and any complex number $z$ with positive 
real part. Herein $\mu_b$ denotes the law of the first passage time to the 
level $b$, given by 
$$\mu_b(dw)={b\over \sqrt{2\pi\, }} w^{-3/2}e\vbox to 9pt{}^{-
{\scriptstyle  b^2\over \scriptstyle 2w}}, $$
on the positive real line.
\bigskip
{\bf Proof of the basic valuation identity:}\quad The proof of the 
basic valuation identity is by reduction to that of [{\bf CJY}, \S\S 4, 5]. 
As a first step, the valuation problem for the Paris option 
of paragraph two is transcribed in in terms of a restarted--at--time--$t$
normalized Brownian motion with drift $W^*$ as follows. With  $B^*$ the 
restarted--at--time--$t$ Brownian motion $B^*(u)=B(t\!+\!u)\!-\!B(t)$, 
consider the Brownian motion with drift:
$$W_u^*={1\over \sigma}\log\Bigl( {S_{t+u}\over S_t}\Bigr) 
=\varpi u\!+\! B_u^*\, ,
\qquad u\in [t,\infty),$$
whose drift coefficient is  
$\varpi=\sigma^{-1}\cdot(r\!-\!\delta\!-\sigma^2/2)$. Recalling 
$b=\sigma^{-1}\cdot\log(L/S_t)$, one has $S_{t+u}<L$ iff $W^*(u)<b$
by construction. Thus define the excursion time $H_b^*=H_{b,D,0}^{*,-}$  
for $W^*$ as follows. If today's value $W^*(0)$ of $W^*$ is 
greater than or equal to $b$ then $H_{b}^*$ is the first time $s$ 
greater than or equal to $D$ such that $W^*(u)<b$ for any time $u$ in the 
interval $(s\!-\!D,s)$ of the positive real line. Suppose 
today's value $W^*(0)$ of $W^*$ is less than $b$. Set $H_{b}^*=d$ if 
$W^*(u)<b$ for any time $u$ in the interval $[0,d)$ remaining  for 
$W^*$ to have stayed a connected period of time of length $D$ below $L$. 
Otherwise, $H_{b}^*$ is the first 
time $s$ greater than or equal to $D$ such that $W^*(u)<b$ for any 
time $u$ in the subinterval $(s\!-\!D,s)$ of the positive real line.
One has  $H_b^*=H_{L,t}\!-\!t$, whence \hbox{$H_{L,t}\le T$ iff
$H_b^*\le \tau=T\!-\!t$}. Using Strong Markov, the valuation problem
for the Paris option thus transcribes  as follows: 
$$ C_{d,i}=e\vbox to 7pt {}^{-r\tau}E^Q\Big[
 \left( S_t\cdot e\vbox to 9pt{}^{\sigma W^*(\tau)}\!-\! K\right)^{+}\cdot 
{\bf 1}_{\{H_b^*\le \tau\}}\Big]. 
$$ 
As a next step, using the particular case [{\bf KS}, p.196f] of the 
Cameron--Martin--Girsanov theorem, change measure from $Q$ to $Q^*$ such
that the drift term in $W^*$ is killed and $W^*$ becomes a $Q^*$--Brownian
motion. Change measure in the expectation expressing the price $C_{d,i}$
and iterate the $Q^*$--expectation thus obtained to get:
$$C_{d,i}=e\vbox to9pt{}^{
        -\bigl( r+{\scriptstyle \varpi^2\over\scriptstyle 2}\bigr)
                         \cdot \tau}
E^*\Big[ {\bf 1}_{\{ H_{b}^*\le \tau\}}
E^*\Big[ e^{\varpi W_\tau^*}\bigl( S_t\cdot e^{\sigma W^*_\tau}\!-\! K\bigr)^+
\Big| \fh _{H_{b}^*}\Bigr]\Bigr].$$
Using $H_b^*$ as $\fh$--stopping time apply the strong Markov 
property of Brownian motion in the conditional expectation. The 
resolvent thus obtained can be further simplified with the random
variables $H_{b}^*$ and $W^*(H_{b}^*)$ being independent. 
With $\mu^*$ the measure for  $W^*(H_{b}^*)$, the basic valuation identity 
follows. 
\goodbreak
\bigskip
{\paragraphfont 7.\quad  An intermediate identity for the proofs of 
the Laplace transforms:} \quad 
Adopt the notation introduced in the previous paragraph six in 
the basic valuation identity for Paris options. Then one has the 
following intermediate identity:
\bigskip
{\rfont Lemma:}\quad{\it For any real number $y$,  one has:
$$ \la\bigl( h_b(\ , y)\bigr)(z)=\int_{\bf R} E^*\Big[ 
e\vbox to 9pt{}^{-z H_{b}^*}\Big] \cdot 
{e\vbox to 9pt{}^{-|x-y|\sqrt{2z\,}}
 \over \sqrt{2z\,}}\,  \mu^*(dx)\, , 
$$
for any complex number $z$ in the right complex half--plane.}
\bigskip
For proving the Lemma interchange the Laplace transform with the two
exponential integrals defining the Paris option density $h_b(u,y)$ 
in paragraph six to get:
$$ \la\bigl( h_b(\ , y)\bigr)(z)= 
\int_{\bf R} 
E^*\bigg[ \int_{H_{b}^*}^\infty 
e^{-zu} {1\over \sqrt{ 2\pi (u\!-\! H_{b}^*)\, }}\cdot 
e\vbox to 9pt{}^{
-{\scriptstyle (x-y)^2\over\scriptstyle 2 (u\!-\! H_{b}^*)}
              }\, du \bigg]\, \mu^*(dx)\, . 
$$
In this last integral successively change variables $w=u\!-\!H_{b}^*$ and 
separate the expectation from the Laplace transform
to get: 
$$\int_{\bf R} 
E^*\Big[  e\vbox to 7pt{}^{ -zH_{b}^*} \Big]\cdot 
  \la \biggl({1\over \sqrt{2\pi u\, }}\cdot 
  e\vbox to 9pt{}^{
-{\scriptstyle (x-y)^2\over\scriptstyle 2 u}
              }\biggr)(z)\,  \mu^*(dx).$$
Using [{\bf Sch}, \S 3] for the Laplace transform completes
the proof of the Lemma. 
\goodbreak
\bigskip
{\paragraphfont 8.\quad Computing the Laplace transform in the case $b$ 
equal to zero:}\quad The computation of the Laplace transforms of paragraph 
six for $h_b$ is by reduction to the case $b$ equal to zero where 
one has to show the {\it key relation\/}:
$$ E^*\Big[ e\vbox to 9pt{}^{
-{\scriptstyle z^2\over\scriptstyle2 }H_0^*}
\Big] \cdot \Psi(z\sqrt{D\, })=\Psi(0)=1,$$
for any complex number $z$ with positive real part. This paragraph 
reviews the key ideas of its proof in [{\bf CJY}, \S 8]. 
The computations in this case $b$ equal to zero are based on the Az\'ema 
martingale and properties of the Brownian meander. 
The {\it Az\'ema martingale\/} is the martingale on $[0,\infty)$ for 
the slow Brownian filtration ${\bf F}^+$ given for any $t\ge 0$ by:
$$\mu_t=(\sgn W_t^*)\sqrt{ t\!-\!g_t\, }.$$
Herein $g_t$ is defined as the supremum over all real numbers $s\le t$ 
such that $W^*(s)=0$. The process given by:
$$m_t(u)={1\over \sqrt{ t\!-\!g_t\, }} \big|
W^*\bigl( g_t+u(t\!-\!g_t)\bigr)\big|, \qquad u\in[0,1]$$
is a {\it Brownian meander\/} and is independent of the $\sigma$--subfield
$\fh^+(g_t)$ of the slow Brownian filtration. Its law is independent of 
$t$ and given by $x\exp(-x^2/2)\cdot {\bf 1}_{(0,\infty)}\, dx$. The 
tautology $W^*(t)=m_t(1)\cdot \mu_t$ thus implies the following identity:
$$
E^*\bigg[ e\vbox to 9pt{}^{ zW^*_t-{\scriptstyle z^2\over\scriptstyle2 }t}
\Big| \fh^{\, +}(g_t)\Big] =       
 e\vbox to 9pt{}^{-{\scriptstyle z^2\over\scriptstyle2 }t}
\cdot \Psi(z\cdot \mu_t)\, , $$
for any positive real number $z$.  Fixing any such $z$, the crucial point 
is to exhibit the product $\psi(-z\cdot \mu(H_0^*))
\exp(-z^2H_0^*/2)$ as a 
martingale. This is achieved by a boundedness argument showing the argument
of the above conditional expectation to be uniformly integrable for any 
time $t$ up to $H_0^*$. Using another 
form of the optional stopping 
theorem it follows that the $Q^*$--expectation of the martingale stopped 
at $H_0^*$ equals the $Q^*$--expectation of it at time zero, whence
$\Psi(0)$ which equals one. With the random variables $H_0$ and 
$W^*(H_0^*)$ independent, the  key relation results. 
This is an identity between functions holomorphic in particular
on the right half--plane. Using the identity theorem, it thus remains 
valid for any complex number $z$ with positive real part. 
\goodbreak
\bigskip   
{\paragraphfont 9.\quad  Computing the Laplace transform for $b$ 
non--positive:}\quad Using \S 7~Lemma the proof of \S 6~Proposition~B 
reduces to compute the expectation of $\exp(-zH_b^*)$, 
for any complex number $z$ with positive real part, and  
the density $\mu^*$.  Following [{\bf CJY}, \S 8.3.3], this is by 
reduction to the case $b$ equal to zero of the previous paragraph.
Indeed, decompose $H_{b}^*$ as follows:
$$ H_{b}^*=T_b+H_0^{**}. $$
Herein $T_b$ is the first passage time of $W^*$ to the level $b$, and 
$H_0^{**}$ is defined as follows. It is the smallest non--negative 
point in time $s$ at which the restarted--at--time--$T_b$ Brownian 
motion $W^{**}(u):=W^*(T_b\!+\!u)\!-\!W^*(T_b)=W^*(T_b\!+\!u)-b$ 
is zero for the first time  after having been less than zero for
a connected period of time of length at least $D$. 
\medskip
To compute the expectation required, notice that $T_b$ and $H_0^{**}$ 
are independent random variables. The conditional expectation 
at time $T_b$ of $\exp(-(w^2/2)H_b^*)$ thus is the product 
of $\exp(-(w^2/2)T_b)$ and the expectation of $\exp(-(w^2/2)H_0^{**})$. 
With the law of $H_0^{**}$ equivalent to that of $H_{0}^*$,
this last expectation is given by the key relation of paragraph eight. 
So it is deterministic in particular. Take expectations to time zero
of this product. The expectation of $\exp(-(w^2/2)T_b)$  that remains
to be computed is 
the Laplace transform of the law of $T_b$ at $w^2/2$. Using  
[{\bf Sch}, \S 3] it is equal to $\exp(-|b|(2\cdot (w^2/2))^{1/2})$ 
and thus to $\exp(+bw)$. With $w=(2z)^{1/2}$ one gets:
$$ E^*\Big[ e\vbox to 7pt{}^{ -z\cdot H_b^*}\Big]
= {1\over \Psi(\sqrt{2Dz\,}\, )}\cdot
E^*\Big[ e\vbox to 7pt{}^{ -\sqrt{2z\, }\cdot T_b}\Big]
= {1\over  \Psi(\sqrt{2Dz\,}\, )}\cdot e\vbox to 7pt{}^{b\sqrt{2z\,}}\, .$$  
To determine $\mu^*$  notice the tautology
$W^*(H_b^*) 
=W^{**}(H_0^{**})\!+\!b$. Thus  
$W^*(H_b^*)\le x$ iff  $W^{**}(H_0^{**})\le x\!-\!b$. From the 
discussion of the  case $b$ equal to zero, the law of 
$W^{**}(H_0^{**})$ is that of the negative of $m_1$. As a consequence:
$$ Q^*\bigl( W^*(H_b^*)\in dx\bigr) =
{\bf 1}_{(-\infty,b]}(x)\cdot (b-x)\cdot e\vbox to 9pt{}^{ -
     {\scriptstyle (x-b)^2\over \scriptstyle 2D}}\cdot {dx\over D}\, .$$
Substitute into  \S 7~Lemma and change variables $w=b\!-\!x$ to
obtain  \S 6~Proposition~B for any complex number 
$z$ with sufficiently big positive real part. Extend to the whole 
right half--plane using the consequence of the key relation of paragraph
eight that $\Psi((2Dz)^{1/2})$ has no zeroes there. This completes the 
proof of \S 6~Proposition~B and the $b$--non--positive part of 
\S 6~Proposition~A. 
\goodbreak
\bigskip
{\paragraphfont 10.\quad   Computing the Laplace transform for $b$ 
positive:}\quad The proof of Proposition~C of paragraph six modifies
the argument of [{\bf CJY}, \S8.3.3]. Using \S 7~Lemma, the first task
is to compute the density $\mu^*$ and the $Q^*$--expectation of 
$\exp(-zH_b^*)$ for any complex number $z$ with positive real part. 
The two cases where the level $b$ is hit before or after time $d$ are 
to be distinguished. Thus decompose the underlying probablity space 
$\Omega$ into the set $A$ on which the first passage time $T_b$ of $W^*$ 
to the level $b$ is less than or equal to $d$ and its 
complement $\Omega\setminus A$ on which $T_b$ is bigger than $d$. 
This induces the decomposition:
$$
E^*\Big[ e\vbox to 7pt{}^{ -zH_b^*}\Big] 
=E^*\Big[ {\bf 1}_A \cdot e\vbox to 7pt{}^{ -zH_b^*}\Big]
+E^*\Big[ {\bf 1}_{\Omega\setminus A} \cdot 
e\vbox to 7pt{}^{ -zH_b^*}\Big]$$
of the $Q^*$--expectation of the random variable $\exp(-zH_b^*)$. 
By construction $H_b^*=d$ on the complement of $A$. So the 
second above summand is $\exp(-zd)Q^*(T_b>d)$. On the set $A$ the level 
$b$ is reached before the critical time $d$ and the clock for the excursion
is reset to zero at the first passage time $T_b$ to the level $b$. 
\hbox{Accordingly one has the decomposition}:
$$ H_{b}^*=T_b+H_0^{**}, $$
with $H_0^{**}$ defined as follows. It is the smallest $s\ge 0$
at which the restarted--at--time--$T_b$ Brownian motion 
$W^{**}(u):=W^*(T_b\!+\!u)\!-\!W^*(T_b)=W^*(T_b\!+\!u)\!-\!b$ 
is zero for the first time  after having been less than zero for
\hbox{a connected period of time of lenght at least $D$}. 
\medskip
With this random variable $H_0^{**}$ independent from $T_b$, the conditional
expectation of the random variable ${\bf 1}_A\cdot\exp(-zH_{b}^*)$ at 
time $T_b$ thus equals $\exp(-zT_b)$ times the expectation of 
$\exp(-zH_0^{**})$.
The key relation of paragraph eight now applies to $H_0^{**}$ and identifies  
this last expectation as  the reciprocal of $\Psi((2zD)^{1/2})$. The 
expectation of ${\bf 1}_A\cdot\exp(-zT_b)$  on the other hand, is given 
by integrating $\exp(-zw)$ from zero to $d$ against the density $\mu_b$ 
of $T_b$. Summarizing, it so follows for any complex number $z$ with 
positive real part:
$$
E^*\big[ e\vbox to 7pt{}^{ -zH_{b}^*}\big] 
= Q^*(T_b>d)e^{-zd}
+ {\int_0^d e^{-zw}\, \mu_b(dw)\over \Psi(\sqrt{2Dz\,}\,)}\, .
$$
For determining $\mu^*$ decompose 
the random variable $W^*(H_{b}^*)$ with respect to $A$:
$$W^*(H_{b}^*)= \bigl(W^**(H_0^{**})\!+\!b\bigr)\cdot{\bf 1}_A
+ W^*(d)\cdot{\bf 1}_{\Omega\setminus A}\, . $$ 
Using the independence of $T_b$ and $H_0^{**}$, the law of the first summand 
is obtained as the convolution of the laws of $T_b$ and $W^{**}(H_0^{**})$, 
whence
$$\eqalign{
\mu_1^*(dx)&=
Q^*\bigl( (W^{**}(H_0^{**})\!+\!b)\in dx; T_b\le d\bigr)\cr 
&=
Q^*(T_b\le d)\cdot {\bf 1}_{(-\infty ,b]}(x)\cdot (b\!-\!x)\cdot
e\vbox to 9pt {}^{ -{\scriptstyle (x-b)^2\over\scriptstyle 2D}  }
\cdot {dx\over D}\, .\cr 
}$$
For the law of the second summand notice $T_b>d$ iff $W^*(t)<b$ for all
$t\le d$, or equivalently, $\max\{ W^*(t)| t\le d\}$ is smaller than $b$. 
Using [{\bf H}, p.9] one so has:
$$\eqalign{
\mu^*_2(dx)&= 
Q^*\bigl( W^*(d)\in dx; T_b>d\bigr)\cr
&= 
e\vbox to 9pt {}^{ -{\scriptstyle x^2\over\scriptstyle 2d}  }
\cdot{dx\over \sqrt{2\pi d\, }}
-
e\vbox to 9pt {}^{ -{\scriptstyle (x-b)^2\over\scriptstyle 2D}  }dx
\cdot{dx\over \sqrt{2\pi d\, }}\, . \cr
}$$
With the expectation factor in \S 7~Lemma seen above to be 
deterministic and independent of the variable $x$, the Laplace transform 
of $h_b(\ ,y)$ at $z$ is given as the product:
$$ 
\la\bigl( h_b(\ , y)\bigr)(z)=
E^*\big[ e\vbox to 7pt{}^{ -zH_{b}^*}\big] 
\int_{\bf R}
{e\vbox to 9pt{}^{-|x-y|\sqrt{2z\,}}
 \over \sqrt{2z\,}}\,  ( \mu^*_1+\mu^*_2)(dx)\, \, . 
$$ 
One is so reduced to compute the two integrals of the second factor 
of this product. For the first of these, on substituting for $\mu^*_1$ and 
changing variables $w=b\!-\!x$, one gets:
$$
\int_{\bf R}{e\vbox to 9pt{}^{-|x-y|\sqrt{2z\,}} \over \sqrt{2z\,}}\,  
\mu^*_1(dx)=
{Q^*(T_b\le d)\over D}\int_0^\infty 
x\cdot {e\vbox to 9pt{}^{-|b-x-y|\sqrt{2z\,}}  
\over \sqrt{2z\,}}\cdot
e\vbox to 9pt{}^{-{\scriptstyle x^2\over\scriptstyle 2D} }
\, dx\, .
$$
For the second of these one analogously obtains:
$$
\int_{\bf R}{e\vbox to 9pt{}^{-|x-y|\sqrt{2z\,}} \over \sqrt{2z\,}}\,  
\mu^*_2(dx)=
{1\over \sqrt{2\pi d\, }}\int_{\bf R}
{e\vbox to 9pt{}^{-|x-y|\sqrt{2z\,}}\over \sqrt{2z\,}}
\Bigl( 
e\vbox to 9pt{}^{-{\scriptstyle x^2\over\scriptstyle 2d} }\
-
e\vbox to 9pt{}^{-{\scriptstyle (x-2b)^2\over\scriptstyle 2d} }\
\Bigr)dx\, .
$$
At this point it remains to identify the two Laplace inverses
of \S 6~Proposition~C to complete its proof. This is based on
the result: 
$$ \la^{-1}\Bigl( {1\over \sqrt{2z\, }}\cdot
e^{-\alpha\sqrt{2z\,}}\Bigr)(u)
={1\over \sqrt{2\pi u\, }}\cdot
e\vbox to 9pt{}^{-{\scriptstyle \alpha^2\over\scriptstyle 2u} }
$$
valid for any positive real numbers $\alpha$ and $u$, and discussed in  
[{\bf Sch}, \S 3]. Applying Fubini's theorem to get the inversion integral 
as inner integral, the improper integral factor of the above $\mu_1^*$ 
integral becomes:
$$\eqalign{
\int_0^\infty 
x
e\vbox to 9pt{}^{-{\scriptstyle x^2\over\scriptstyle 2D} }
\cdot \la^{-1}\biggl(e\vbox to 6pt{}^{-zd}
{e\vbox to 9pt{}^{-|b-x-y|\sqrt{2z\,} }\over\sqrt{2z\,}}
\biggr)\! (u)\, dx\, . 
\cr }
$$
Successively apply the shifting theorem for the Laplace transform 
of paragraph three to take care of the factor $\exp(-zd)$ and the above
Laplace inversion formula with $\alpha$ equal to $|b\!-\!x\!-\!y|$ to 
arrive at: 
$$
{1\over D}\cdot {\bf 1}_{(d,\infty)}(u){1\over \sqrt{ 2\pi(u\!-\!d)\, }}
\int_0^\infty  
x\cdot e\vbox to 9pt{}^{-\bigl({\scriptstyle x^2\over \scriptstyle 2D}+
    {\scriptstyle (b-x-y)^2\over\scriptstyle 2(u\!-\!d)}\bigr)}
dx\, .
$$ 
This integral is seen to be the value of the function 
$h_{b,3}$ at $u$ and $y$. 
Analogously, the  Laplace 
inverse at any positive real number $u$ of the integral factor 
of the above $\mu^*_2$ integral is: 
$$
{1\over \sqrt{2\pi d\, }} {\bf 1}_{(d,\infty)}(u)
{1\over \sqrt{ 2\pi(u\!-\!d)\, }}
\int_{\bf R} 
 e\vbox to 9pt{}^{ -{\scriptstyle (x-y)^2\over\scriptstyle 2(u\!-\!d)}  }
\Bigl( e\vbox to 9pt{}^{-{\scriptstyle x^2\over \scriptstyle 2d}}-
   e\vbox to 9pt{}^{-{\scriptstyle (x-2b)^2\over\scriptstyle 2d}}
\Bigr)dx\,
$$
and is seen to be equal to the value of the function $h_{b,4}$ at 
$u$ and $y$. This completes the proof of \S 6~Proposition~C.
\goodbreak  
\vskip.7cm
\centerline{\headingfont Part III \qquad Analytic properties of the function
$ \Psi$}   
\vskip.3cm
{\paragraphfont 11.\quad Analytic properties of $\Psi$:}\quad This 
paragraph collects and discusses pertinent analytical properties
of the function $\Psi$ for any complex number $w$ recalled to be given by:
$$ \Psi(w)=\int_0^\infty x\cdot e\vbox to 9pt{}^{ 
-{\scriptstyle x^2\over \scriptstyle2}+wx}\, dx\, .$$
Developing the linear exponential factor
of the integrand in its series and integrating the resulting series 
term by term, gives the following {\it series expansion}:
$$ \Psi(w)=\sum\nolimits_{n=0}^\infty
{ 2^{\scriptstyle n\over \scriptstyle2}\over n!}
\cdot \Gamma\Bigl({n\!+\!2\over 2}\Bigr)\cdot w^n \, . $$
This series is absolutely convergent for any complex number $w$, and 
convergence is uniform on compact sets. Rearrange it in its even and
odd order terms. Replacing $w$ by its negative leaves unchanged the 
even part and produces a minus sign in the odd part. Using the duplication
identity for the gamma function and redeveloping the square order 
exponential series that results as a factor, then gives the following 
{\it key identity\/}:
$$\Psi(w)=\Psi(-w)+\sqrt{2\pi \, }\cdot w\cdot\
e\vbox to 9pt{}^{{\scriptstyle w^2\over \scriptstyle 2}}$$ 
connecting the values of $\Psi$ on the right half--plane with those 
on the left half--plane. 
\goodbreak
\bigskip
{\paragraphfont 12.\quad A uniform asymptotic expansion:}\quad  The aim 
of this paragraph is to prove for the function $\Psi$ given by 
$\Psi(w)=\int_{\raise1.5pt\hbox{$\scriptstyle 0$}}^\infty 
x\exp(-x^2/2+wx)\, dx$, 
for any complex number $w$, the following uniform asymptotic expansion.
\bigskip
{\rfont Lemma:}\quad {\it For any complex number $w$ with positive 
real part, one has the expansion:
$$\Psi(-w)={1\over \sqrt{\pi\,}} \sum_{k=1}^N (-1)^{k+1} 
\cdot 2^k\cdot \Gamma\Bigl( k\!+\!{1\over 2}\Bigr) \cdot {1\over w^{2k}}
+R_{N+1}(w),$$
for any positive integer $N$, whose remainder term $R_{N+1}$ 
satisfies the estimate:}       
$$\big|R_{N+1}(w)\big|\le 
2\cdot {(2N\!+\!1)!\over N!}\cdot {1\over |w|^{2(N+1)}}\, .$$
\medskip
{\rfont Corollary:}\quad{\it For any complex number $z$ not a non--positive 
real number, one has:
$$\Psi(-\sqrt{z\, }\, )={2\over z}+R_2(\sqrt{z\, }\, )$$
with $R_2(z^{1/2})$ being big Oh in $|z|^{-2}$ for the absolute value 
of $z$ going to infinity.}
\bigskip
The Corollary is an immediate consequence of the Lemma. Recall that
the principal branch of the logarithm has been chosen.  
The proof of the Lemma is in two steps. 
As a first step, establish the asymptotic expansion with the remainder term
satisfiying the estimate:
$$\big|R_{N+1}(w)\big|\le
{C_N\over \cos^{2(N+1)}(\pi/2\!-\!\delta_z)}\cdot {1\over |w|^{2(N+1)}}
\qquad\hbox{where}\qquad
C_N={(2N\!+\!1)!\over N!\cdot 2^N}\, , $$
if $|\arg(w)|\le \pi/2\!-\!\delta_w$ for $\delta_w>0$. For this develop 
the square exponential factor of the integrand for $\Psi(-w)$ into its
Taylor series up to the $(N\!\!-\!\!1)$--st term and integrate the resulting 
sum term by term. Using the duplication formula for the gamma function
the coefficients of the asymptotic expansion follow on shifting the 
summation index by one. Proceeding analogously, the absolute value of 
the remainder term can be majorized by $C_N$ times the reciprocal of the 
$2(N\!\!+\!\!1)$--st power of the real part of $w$. With $w$ in the right
half--plane, $\re(w)$ equals $|w|\cdot \cos(\theta)$ for an angle $\theta$
of absolute value less than or equal to $\pi/2\!-\!\delta_w$ as above. 
Thus $\cos (\theta)$ is bigger than the positive constant 
$\cos(\pi/2\!-\!\delta_w)$ completing the proof of the first step. 
\medskip
As a second step, extend the above asymptotic expansion by 
showing the above estimate to be valid also for any complex numbers $w$
such that $\pi/4<|\arg(w)|\le 3\pi/4\!-\!\delta_w$.  This is based on 
the following integral representation:
$$\Psi(-w)=\int_0^{e^{i\theta}\cdot \infty} \xi\cdot e\vbox to 9pt{}^{
 -{\scriptstyle \xi^2\over \scriptstyle 2}-w\xi }\, d\xi\, , $$
for any complex number $w$. Herein $\theta$ is any fixed any angle 
of absolute value less than $\pi/4$ and the integral is over the path 
sending any positive real number $R$ to $R\cdot \exp(i\theta)$. 
This representation is proved using Cauchy's theorem. Indeed, the absolute 
value of the above integrand is majorized by 
$|\xi|\cdot \exp( -\re(\xi^2)/2-\re(w\xi))$. It thus converges uniformly 
to zero with the radius $R$ going to infinity on that part of the circle 
in the origin of radius $R$ which is parametrized by angles of absolute 
value less than $\pi/4$. 
\medskip
Now fix any positive real number $\varepsilon$ smaller than $\delta_w$. 
If $\arg(w)$ is negative, choose the angle $\theta=+\pi/4\!-\!\varepsilon$ 
in the above integral representation. Otherwise let 
$\theta=-\pi/4\!+\!\varepsilon$.  
In the sequel the case $\arg(w)$ positive is considered. In the above 
integral representation then develop the square exponential factor of the   
integrand into its Taylor series up to the $(N\!-\!1)$--st order.
In the identity:
$$\int_0^{e^{i\theta}\cdot \infty}\xi^{2k+1}e^{-w\xi}\, d\xi=
e^{2(k+1)i\theta }\int_0^\infty R^{2k+1}e^{-we^{i\theta}R}\, dR
$$
write $w\exp(i\theta)=|w|\exp(i(\theta \!+\!\arg(w)))$. Notice that
the angle $\theta \!+\!\arg(w))$ is between $\theta \!+\!\pi/4=\varepsilon$ 
and 
$\theta\!+\!3\pi/4\!-\!\delta_w=\pi/2\!-\!(\delta_w\!-\!\varepsilon)$. Thus 
it is smaller than $\pi/2$ in particular. As a consequence, the 
real part of $\eta=w\exp(i\theta)$ is positive. Recalling  that 
$\Gamma(a\!+\!1)$ is equal to the integral 
$\eta^{a+1}\int_{\raise1.5pt\hbox{$\scriptstyle 0$}}^\infty 
\exp(-\eta t)t^a\, dt$, the above integral is thus seen 
to be equal to 
$w^{\raise-1pt\hbox{$\scriptstyle -2(k+1)$}}\cdot \Gamma(2(k\!+\!1))$.  
The coefficients of the asymptotic expansion follow. 
Proceeding analogously, the 
absolute value of the remainder term $R_{N+1}(w)$ is majorized by $C_N$ 
times the reciprocal of the $2(N\!+\!1)$--st power of 
$|w|\cdot \cos(\theta\!+\!\arg(z))$. Recalling the above bounds for its
argument, the cosine herein is positive and bounded from below by 
$\cos(\pi/2\!-\!(\delta_w\!-\!\varepsilon))$, for any positive $\varepsilon$
less than $\delta_w$. Thus it is bounded from below by 
$\cos(\pi/2\!-\!\delta_w)$, completing the proof of the second step. 
\medskip
To complete the proof of the uniform asymptotic expansion, apply the 
estimate of the first step with any $\delta_z$ strictly between $\pi/4$ 
and $\pi/2$ if $|\arg(w)|$ is less than $\pi/4$. Otherwise, apply the 
estimate of the second step with $\delta_w$ strictly between  $\pi/4$ 
and $\pi/2$. In both cases, $\cos(\pi/2\!-\!\delta_w)$ is minorized by 
$\cos(\pi/4)=2^{-1/2}$. Cancelling the powers of two one gets the factor
two. This completes the proof of the uniform asymptotic expansion.
\goodbreak
\bigskip
\medskip
{\paragraphfont References:}
\medskip
{\ninerm 
\baselineskip9.75pt
\parein{35pt}{\ninebf [CJY]} M. Chesney, M. Jeanblanc--Picqu\'e, M. Yor, 
Brownian excursions and Parisian barrier options, Adv. Appl. Probability
{\ninebf 29}(1997), 165--184
\paraus 
\parein{35pt}{\ninebf [D]} G. Doetsch, {\nineit Handbuch der Laplace 
Transformation\/} I, Birkh\"auser Verlag 1971
\paraus
\parein{35pt}{\ninebf [H]} J.M. Harrison,  {\nineit Brownian motion and
stochastic flow systems\/}, Krieger reprint 1990
\paraus
\parein{35pt}{\ninebf [KS]} I. Karatzas, S.E. Shreve, {\nineit Brownian motion
and stochastic calculus\/}, 2nd ed., GTM 113, Springer 1991
\paraus
\parein{35pt}{\ninebf [L]} N.N. Lebedev,  {\nineit Special functions and
their applications\/}, Dover Publications 1972
\paraus
\parein{35pt}{\ninebf [Sch]}  M. Schr\"oder, Fact and fantasy in recent
applications of transform methods to the valuation of exotic options, 
Fakult\"at f\"ur Mathematik und Informatik, Universit\"at Mannheim,  
September 1999
\paraus
\parein{35pt}{\ninebf [SP]}  M. Schr\"oder, On the valuation of 
Paris options: the first standard case, Fakult\"at f\"ur Mathematik 
und Informatik, Universit\"at Mannheim,  November 1999
\paraus
}
\vfill
\eject
\headline={  }
\hoffset-.5cm
\nopagenumbers
\baselineskip20pt
\centerline{  }
\vskip4truecm
\centerline{\gr  On the valuation of Paris barrier options:}
\vskip2pt
\centerline{\gr foundational results} 
\vskip.5cm
\baselineskip14pt
\centerline{\gross Michael Schr\"oder}
\centerline{\grossrm (Mannheim)}
\vskip.3cm
\centerline{\grossrm November 1999}
\bye